\documentclass{amsart}

\newcommand{\bt}{\begin{Theorem}}
\newcommand{\et}{\end{Theorem}}
\newcommand{\bi}{\begin{itemize}}
\newcommand{\ei}{\end{itemize}}
\newcommand{\bea}{\begin{eqnarray}}
\newcommand{\ba}{\begin{array}}
\newcommand{\eea}{\end{eqnarray}}
\newcommand{\ea}{\end{array}}

\newcommand{\f}{\frac}
\newcommand{\s}{\sigma}

\newcommand{\Si}{\Sigma}
\newcommand{\iy}{\infty}
\newcommand{\what}{\widehat}
\newcommand{\wtilde}{\widetilde}

\newtheorem{Definition}{Definition}[section]
\newtheorem{Theorem}[Definition]{Theorem}
\newtheorem{Lemma}[Definition]{Lemma}

\newtheorem{Corollary}[Definition]{Corollary}

\newtheorem{Remark}[Definition]{Remark}

\newcommand{\be}{\begin{equation}}
\newcommand{\ee}{\end{equation}}

\newcommand{\newsection}{\setcounter{equation}{0}}

\newcommand{\R}{\mathbb R}%
\newcommand{\C}{\mathbb C}%
\newcommand{\Z}{\mathbb Z}%
\newcommand{\ff}{\mathfrak f}%

\begin{document}
\baselineskip24pt
\author[Sarkar]{Rudra P. Sarkar}
\address[R. P. Sarkar]{Stat-Math Unit, Indian Statistical
Institute, 203 B. T. Rd., Calcutta 700108, India, E-mail:
rudra@isical.ac.in}
\author[Sengupta]{Jyoti Sengupta}
\address[J. Sengupta]{School of Mathematics, Tata Institute of Fundamental
Research, Homi Bhabha Rd.,  Mumbai 400005, India, E-mail:
sengupta@math.tifr.res.in}

\title{Beurling's Theorem for $SL(2,\R)$}
\begin{abstract}
We prove  Beurling's theorem for the full group $SL(2,\R)$. This
is the {\em master theorem} in the quantitative uncertainty
principle as all the other theorems of this genre follow from it.
\end{abstract}
\maketitle
\section{Introduction} Our starting point is the following theorem
of H\"ormander (\cite{Hor}):
\begin{Theorem} {\em (H\"ormander 1991)} Let $f\in L^1(\R)$. If
$$\int_\R\int_\R|f(x)\wtilde{f}(y)|e^{|xy|}dx\,dy<\infty$$
where $\wtilde{f}(\lambda)=\int_\R f(t)e^{-i\lambda t}dt$, then
$f=0$ almost everywhere. \label{thm-hormander}
\end{Theorem}

H\"ormander attributes this theorem to A. Beurling. We will follow
his practice and call theorem 1.1 as {\em Beurling's theorem}.
This theorem is an instance of the quantitative uncertainty
principle (QUP) (see \cite{FS}). We recall that the uncertainty
principle is the paradigm in harmonic analysis which says that a
function and its Fourier transform cannot both be very rapidly
decreasing at infinity. Some other well known QUP's like
Cowling-Price theorem, Gelfand-Shilov theorem, Morgan's theorem
and Hardy's theorem (see section 4 for the precise statements)
become corollaries of this theorem. Thus Beurling's theorem can be
regarded as the ``Master Theorem'' in the context of the
uncertainty principle.

In recent years mathematicians have taken up the QUP problems on
semisimple Lie groups and on Riemannian symmetric spaces and
proved versions of Hardy's, Cowling-Price, Gelfand-Shilov and
Morgan's theorems (see \cite{CSS, FS, NR, Sa1, Se, SiSu, T1} and
the references therein). However, Beurling's theorem has not yet
been proved for any semisimple Lie group. The aim of this article
is to prove Beurling's theorem for the full group $SL(2, \R)$. The
other four theorems  mentioned above will follow from it.

Our principal tool is the Abel transform. This is a   well-known
device for handling  $K$-biinvariant functions on $G$. In this
article we extend its scope to bi $K$-finite functions  and use
this in a crucial way to obtain our result (see Theorem 3.1) where
no $K$-finiteness restriction is imposed.

It may be noted that this is the first instance when the discrete
series representations appear in the hypothesis in a theorem of
uncertainty. See Remark 4.1 for further details.

 The plan of the article is as follows. In the next section we
record some preliminary material and set up the necessary
notation.  In section 3 we give the statement and proof of our
main result. In section 4 we indicate how the other QUP's viz
Hardy's theorem, Morgan's theorem etc. follow from our main
result.


\section{notation and preliminaries}\newsection
The letter $C$ will denote a positive constant, not necessarily
the same at each occurrence. We will mainly use the notation of
\cite{Ba} with a few exceptions which we will mention here. For
ready reference we will also quote from \cite{Ba} the things which
we will frequently use.
Let $G=SL(2,\R)$.  Let
$$k_\theta= \left (\begin{array}{cc}\cos\theta &\sin\theta\\
-\sin\theta &\cos\theta\end{array}\right ),
a_t=\left(\begin{array}{cc}e^t&0\\0&e^{-t}\end{array}\right )
 \mbox{ and }
n_\xi=\left(\begin{array}{cc}1&\xi\\0&1\end{array}\right ).$$ Then
$K=\{k_\theta\, |\, \theta\in [0, 2\pi)\}$, $A=\{a_t\,|\, t\in
\R\}$ and $N=\{n_\xi\, |\, \xi\in \R \}$ are three particular
subgroups of $G$ of which $K$ is a maximal compact subgroup
$SO(2)$ of $G$. It is clear from the above that both $A$  and $N$
can be identified with $\R$. Let $G=KAN$ be an Iwasawa
decomposition of $G$ and for $x\in G$, let $x=k_\theta a_t n_\xi$
be its corresponding decomposition. We will write $H(x)$ for $t$
and $K(x)$ for $k_\theta$. Clearly $H$ is left $K$-invariant and
right $N$-invariant. The Haar measure $dx$ of $G$ splits according
to this decomposition as $dx=e^{2t}\,dk\,dt\,dn$ where
$dk=dk_\theta=\frac{d\theta}{2\pi}$ is the normalised Haar measure
of $K$ and $dn=dn_\xi=d\xi$ as well as $da=da_t=dt$ are  both
Lebesgue measures on $\R$.

 We also recall that $G$ has Cartan
decomposition $ G=K\overline{A^+}K$, $x=k_1a_tk_2$ where $k_1,k_2
\in K, t\geq 0$. The Haar measure of $G$ splits according to this
decomposition as $dx = dk_1\, \sinh 2t \,dt\, dk_2$. Let
$\sigma(x)=\sigma(k_1a_tk_2)=|t|$. In fact $\sigma(x)=d(xK, o)$,
where $o=eK$ is the `origin' of the symmetric space $G/K$ and $d$
is the distance function on $G/K$.

Let $\what{K}=\{e_n|n\in \Z\}$ be the set of continuous characters
of $K$, where $e_n(k_\theta)=e^{in\theta}$. Instead of $e_n$, by
abuse of language, we will call the integers $n$ as $K$-types. A
complex valued function $f$ on $G$ is said to be of left
(respectively right) $K$-type $n$ if $f(kx)= e_n(k)f(x)$
(respectively $f(xk)= e_n(k)f(x)$) for all $k\in K$ and $x\in G$.
A function is   of type $(m,n)$ if its left $K$-type is $m$ and
right $K$-type is $n$.  For a suitable function $f$, the
$(m,n)$-th isotypical component of $f$ is denoted by $f_{m,n}$ and
this is given by:
\begin{equation}
\int_K\int_K\overline{e_m(k_1)}\,\,\overline{e_n(k_2)}f(k_1xk_2)dk_1\,dk_2
=\int_K\int_Ke^{-im\theta}e^{-in\phi}f(k_\theta
xk_\phi)dk_\theta\, dk_\phi. \label{projection-mn}
\end{equation}
It can be verified that $f_{m,n}$ is itself a function of type
$(m,n)$ and $f_{m,n}\equiv 0$ when $m$ and $n$ are of opposite
parity. The function $f$ can be decomposed as $f=\sum_{m,n\in
\Z}f_{m,n}$. In fact when $f\in C^\infty(G)$ this is an absolutely
convergent series in the $C^\infty$-topology. When $f\in L^p(G)$,
$p\in [1, \infty)$, the equality is  in the sense of distribution.

Let $\mathfrak a$ be the Lie algebra of $A$. Let $\mathfrak a^*$
be the real dual of $\mathfrak a$ and $\mathfrak a^*_\C$ be the
complexification of $\mathfrak a^*$. Then $\mathfrak a^*$ and
$\mathfrak a^*_\C$ can be identified with $\R$ and $\C$
respectively via $\rho$, the half-sum of the positive roots, i.e.
$ \rho=1$ under this identification. Let $M$ be $\{\pm I\}$, where
$I$ is the $2\times 2$ identity matrix. The unitary dual of $M$ is
$\what{M}=\{\sigma_+,\sigma_-\}$ of which $\sigma_+$ is the
trivial representation of $M$ and $\sigma_-$ is the only
nontrivial unitary irreducible representation of $M$. Let
$\Z^{\sigma^+}$ (respectively $\Z^{\sigma^-}$) be the set of even
(respectively odd) integers.

For $\sigma\in \what{M}$ and $\lambda\in \mathfrak a^*_\C=\C$, let
$(\pi_{\sigma,\lambda},H_\sigma)$ be the principal series
representation of $G$ where $H_\sigma$ is the subspace of $L^2(K)$
generated by the orthonormal set $\{e_n|n\in \Z^\sigma\}$. The
representation $\pi_{\sigma,\lambda}$ is normalized so that it is
unitary if and only if $\lambda\in i\mathfrak a^*= i\R$. In fact
(\cite[4.1]{Ba}):
\begin{equation}
(\pi_{\sigma,
\lambda}(x)e_n)(k)=e^{-(\lambda+1)H(x^{-1}k^{-1})}e_{-n}(K(x^{-1}k^{-1})).
\label{principal-series-action}
\end{equation} For every $k\in \Z^*$, the set of nonzero integers,
there is a discrete series representation $\pi_k$ which occurs as
a subrepresentation of $\pi_{\sigma, |k|}$ so that $k\in
\Z\setminus\Z^{\sigma}$.  For $m,n\in \Z^\sigma$ and $k\in
\Z\setminus\Z^{\sigma}$, let $\Phi^{m,n}_{\sigma, \lambda}(x)=
\langle \pi_{\sigma,\lambda}(x)e_m, e_n\rangle$ and
$\Psi^{m,n}_{k}(x)= \langle \pi_{k}(x)e^k_m, e^k_n\rangle_k$, be
the matrix coefficients of the principal series and discrete
series representations respectively, where $\{e^k_n\}$ are the
renormalised basis and $\langle \ ,\  \rangle_k$ is the
renormalised inner product of $\pi_k$ (see \cite[p. ~20]{Ba}). In
particular $\Phi_{\sigma_+, \lambda}^{0,0}$ is clearly the
elementary spherical function, which we also denote by
$\phi_\lambda$. For details of the parametrization of the
representations $\pi_{\sigma,\lambda}$ and $\pi_k$ and their
realizations we refer to \cite{Ba}.

For a function $f \in L^1 (G)$, let $\what{f}(\sigma, \lambda)$
and $\what{f}(k)$ denote its (operator valued) principal and
discrete Fourier transforms at $\pi_{\sigma, \lambda}$ and $\pi_k$
respectively. Precisely:
$$\what{f}(\sigma, \lambda)=\int_Gf(x)\pi_{\sigma,\lambda}(x^{-1})dx
\mbox{ and } \what{f}(k)=\int_Gf(x)\pi_{k}(x^{-1})dx.$$ The
$(m,n)$-th matrix entries of $\what{f}(\sigma, \lambda)$ and
$\what{f}(k)$ are denoted by $\what{f}(\sigma,\lambda)_{m,n}$ and
$\what{f}(k)_{m,n}$ respectively. Thus $\what{f}(\sigma,
\lambda)_{m,n}=\langle \what{f}(\sigma, \lambda)e_m,
e_n\rangle=\int_G f(x)\Phi^{m,n}_{\sigma,\lambda}(x^{-1})dx$ and
$\what{f}(k)_{m,n}=\int_G f(x)\Psi^{m,n}_k(x^{-1})dx$. As $\int_G
f(x)\Phi^{m,n}_{\sigma,\lambda}(x^{-1})dx=\int_G
f_{m,n}(x)\Phi^{m,n}_{\sigma,\lambda}(x^{-1})dx$,  clearly,
$\what{f}(\sigma, \lambda)_{m,n} = \what{f}_{m,n}(\sigma,
\lambda)$. Similarly $\what{f}(k)_{m,n} = \what{f}_{m,n}(k)$.
Henceforth we will not distinguish between $\what{f}(\sigma,
\lambda)_{m,n}$ (respectively $\what{f}(k)_{m,n})$ and
$\what{f}_{m,n}(\sigma, \lambda)$ (respectively
$\what{f}_{m,n}(k)$). Notice that integers $m,n$ of the same
parity uniquely determine a $\sigma\in \what{M}$ by $m,n\in
\Z^\sigma$. Therefore we may sometimes omit the obvious $\sigma$
and write $\Phi_\lambda^{m,n}$ for $\Phi_{\sigma, \lambda}^{m,n}$
and $\what{f}_{m,n}(\lambda)$ for $\what{f}_{m,n}(\sigma,
\lambda)$.

From  (\ref{principal-series-action}) it follows that for
$\lambda\in \C$, $\sigma\in\what{M}$ and $m,n\in \Z^\sigma$,
\begin{equation}\Phi^{m,n}_{\sigma, \lambda}(x)=
\int_Ke^{-(\lambda+1)H(x^{-1}k^{-1})}e_{-m}(K(x^{-1}k^{-1}))e_n(k^{-1})dk.
\label{expression-Phi}\end{equation} Hence,
$|\Phi^{m,n}_\lambda(x)|\le \int_Ke^{-(\Re
\lambda+1)H((x^{-1}k^{-1})}dk=\Phi^{0,0}_{\Re \lambda}=\phi_{\Re
\lambda}$, where $\Re \lambda$ stands for the {\em real part} of
$\lambda$. It is well known (\cite{He2}) that
$|\phi_{\lambda}(x)|\le 1$ for $x\in G$ and $\lambda\in \C$ with
$|\Re\lambda|\le 1$. Combining this with the following two
estimates (\cite[proposition 4.6.4]{GV}  and \cite[3.2]{Ba}):
\begin{equation}
|\phi_{\lambda}(x)|\le e^{|\Re \lambda|\sigma(x)}\Xi(x) \mbox{ for
} \lambda\in \C\mbox{ and } \Xi(x)\le C (1+\sigma(x))
e^{-\sigma(x)} \label{basic-estimates}
\end{equation}
where $\Xi=\Phi^{0,0}_0=\phi_0$,  we get
$|\Phi^{m,n}_\lambda(x)|\le C e^{|\Re
\lambda|\sigma(x)}(1+\sigma(x)) e^{-\sigma(x)}\le  C e^{|\Re
\lambda|\sigma(x)}$, $\forall \lambda\in \C$.

It follows from the estimates above and Morera's theorem that
$\what{f}_{m,n}(\lambda)$ is a holomorphic function in the
Helgason-Johnson strip, $\{\lambda\in \C\mid |{\Re \lambda}|<1\}$.
In particular the restriction of $\what{f}_{m, n}$ to the
imaginary axis is a (complex valued) real analytic function and
hence its zeros form a set of Lebesgue measure zero.
\vspace{.07in}

We conclude this section with a brief discussion of  Plancherel
measure (see \cite{L}). This measure on the unitary principal
series representations (parametrized by $i\R$) is
$d\mu(\sigma,\lambda)=\mu(\sigma, \lambda)d\lambda$ where
\begin{equation}
\ba{lllclll}\mu(\sigma_+, i\xi)&=&(\frac {\xi}
{4\pi})\tanh(\frac{\xi\pi}2)&\mbox{ and }\,\mu(\sigma_-,
i\xi)&=&(\frac {\xi} {4\pi})\coth(\frac{\xi\pi}2)\mbox{ for }
\xi\in \R.\ea \label{plancherelmeasure}
\end{equation}  Here again we may omit $\sigma$ and write
$\mu(\lambda)$ for $\mu(\sigma, \lambda)$, when there is no
confusion about $\sigma\in \what{M}$. The Plancherel measure on
the discrete series is given by $\mu(\pi_k)=\frac{|k|}{2\pi}$, for
$k\in \Z^*$.

\section{Statement and Proof of the main theorem }\newsection
In the proofs of the theorems, Lemmas etc. we will use Fubini's
theorem freely without explicitly mentioning it every time.
\begin{Theorem}
Let $f\in L^2(G)$. If
\begin{equation}\int_G\int_{i\R}|f(x)|\|\what{f}(\sigma,
\lambda)\|_2\phi_{|\lambda|}(x)dxd\mu(\sigma, \lambda)< \infty.
\label{beurling-condition-0-main}
\end{equation}
for all $\sigma\in \what{M}$ and
\begin{equation}
\sum_{k\in \Z^*}\f{|k|}{2\pi}\int_G
|f(x)|\|\what{f}(k)\|_2\phi_{|k|}(x)dx  <\infty
\label{beurling-discrete-condition-main}
\end{equation}
then $f=0$ almost everywhere. Here $\|\ \cdot \ \|_2$ is the
Hilbert-Schmidt norm. \label{thm-sl2-1}
\end{Theorem}

As $\phi_{|\lambda|}(x)\le e^{|\lambda|\sigma(x)}\Xi(x)$ (see
(\ref{basic-estimates})), we have the following immediate
corollary:
\begin{Corollary}
Let $f\in  L^2(G)$. Suppose
\begin{equation}\int_G\int_{i\R}|f(x)|\|\what{f}(\sigma,
\lambda)\|_2e^{|\lambda|\sigma(x)}\Xi(x)dxd\mu(\sigma, \lambda)<
\infty. \label{beurling-condition-0-main-reform1}
\end{equation} for all $\sigma\in \what{M}$
and
\begin{equation} \sum_{k\in \Z^*}\f{|k|}{2\pi}\int_G
|f(x)|\|\what{f}(k)\|_2e^{|k|\sigma(x)}\Xi(x)dx <\infty.
\label{beurling-discrete-condition-main-reform1}
\end{equation}
Then $f=0$ almost everywhere. \label{thm-reform1}
\end{Corollary}

\begin{Remark}{\em
Note that a naive analogue of Beurling's theorem would use the
weight $e^{|\lambda|\sigma(x)}$. Instead here we use
$\phi_{|\lambda|}(x)$ which has less decay (since
$\phi_{|\lambda|}(x)\le e^{|\lambda|\sigma(x)}\Xi(x)$ and
$\Xi(x)\le 1$) and hence obtain a stronger result. We also feel
that the  formulation of the theorem is  natural as
$\phi_\lambda(x)$ plays the role of $e^{i\lambda x}$, at least for
the $K$-biinvariant functions on $G$.}

\end{Remark}
 The basic strategy of our proof of theorem \ref{thm-sl2-1} is to
reduce the theorem to the Euclidean situation by using the Abel
transform. We shall therefore begin with a short discussion on the
Abel transform.
 For $f\in L^1(G)$ of type $(m,n)$, we define the Abel
transform
$${\mathcal A}f(t)=e^{t}\int_Nf(a_tn)dn.$$ Therefore
\begin{equation}|{\mathcal A}f(t)|\le
e^{t}\int_N|f(a_tn)|dn={\mathcal A}|f|(t) \label{mod-f-abel}
\end{equation}
Since $|f|$ is integrable and $K$-biinvariant we have ${\mathcal
A}|f| \in L^1({\R})$. Furthermore ${\mathcal A}|f|$ is an even
function of $t$ ( vide \cite{GV}). Therefore ${\mathcal A}f\in L^1
({\R})$. We need the following lemma:
\begin{Lemma}
Let $\sigma\in \what{M}$ and let $f\in L^1(G)_{m,n}$ for some
$m,n\in \Z^\sigma$. Then
$\what{f}_{m,n}(\sigma,\lambda)=\wtilde{{\mathcal A}f}(-i\lambda)$
for $\lambda\in i\R$. Here  $\wtilde{{\mathcal A}f}(\nu)=\int_\R
{\mathcal A}f(t)e^{-i\nu t}dt$. \label{abel-triangle}
\end{Lemma}
\begin{proof} For reason mentioned in section 2 we will  omit
$\sigma$ and write $\Phi^{m,n}_\lambda$ for $\Phi^{m,n}_{\sigma,
\lambda}$ and $\what{f}(\lambda)$ for $\what{f}_{m,n}(\sigma,
\lambda)$. From (\ref{expression-Phi}) we have for $\lambda\in i
\R$,
$$\ba{lll}\what{f}_{m,n}(\lambda)&=&\int_Gf(x)\int_Ke^{-(\lambda+1)H(xk^{-1})}
e_{-m}(K(xk^{-1}))e_n(k^{-1})dk\,dx\\ \\
&=&\int_K\int_Gf(x)e^{-(\lambda+1)H(xk^{-1})}
e_{-m}(K(xk^{-1}))e_n(k^{-1})dx\,dk.\ea$$

Substituting $k^{-1}yk$ for $x$, we get,\\
$\what{f}_{m,n}(\lambda)=\int_K\int_Gf(k^{-1}yk)e^{-(\lambda+1)H(k^{-1}y)}
e_{-m}(K(k^{-1}y))e_n(k^{-1})dy\,dk$ as the Haar measure of $G$ is
invariant under this substitution.

Since  $f(k^{-1}yk)=f(y)e_m(k^{-1})e_n(k)$,
$e_{-m}(K(k^{-1}y))=e_m(k)e_{-m}(K(y))$ and $H(k^{-1}y)=H(y)$, we
have $\what{f}_{m,n}(\lambda)=\int_Gf(y)e^{-(\lambda+1)H(y)}
e_{-m}(K(y))dy.$

Using the Iwasawa decomposition $G=KAN$ and the identification of
$A$ and $\R$, we obtain,
$$\ba{lll}\what{f}_{m,n}(\lambda)&=&\int_K\int_\R\int_Nf(ka_tn)e^{-(\lambda+1)t}e_{-m}(k)dk\,
e^{2t}\,dt\,dn\\ \\&=&\int_\R\int_Nf(a_tn)e^{-(\lambda+1)t}e^{2t}\,dt\,dn\\
\\
&=&\int_\R e^t\int_Nf(a_tn)dn e^{-\lambda t}dt\\ \\
&=&\int_\R{\mathcal A}(f)(t)e^{-i(-i\lambda) t}dt.\ea$$ Thus
$\what{f}_{m,n}(\lambda)=\wtilde{{\mathcal A}f}(-i\lambda).$
\end{proof}

Note that the lemma above is valid,  for any $\lambda\in \C$ for
which both sides of the equality are well-defined.

Looking back at theorem \ref{thm-hormander} we see that it can be
rewritten as: For $g\in L^1(\R)$, if
\begin{equation}\int_\R
M(g)(\lambda)|\wtilde{g}(\lambda)|d\lambda<\infty,
\label{modified-hormander}
\end{equation} where
$M(g)(\lambda)=\int_\R|g(x)|e^{|\lambda||x|}dx$ then $g=0$ almost
everywhere.


 With this preparation we
are now ready to prove theorem 3.1
\begin{proof} We shall divide the proof in a few steps for
convenience. Before proving Step 1, let us note that if
$\what{f}(\sigma, \cdot)\equiv 0$ on $i\R$, for all $\sigma\in
\what{M}$, then the Fourier transform of $f$ (hence of $f_{m,n}$)
is supported on the discrete series representations. In this case
we can directly go to Step 4.

\noindent{\bf Step 1:} Let us fix a $\sigma\in \what{M}$ such that
$\what{f}(\sigma, \cdot)\not\equiv 0$ on $i\R$. (If there is no
such $\sigma$ then we go to Step 4 as mentioned above.) In this
step we shall show that $f\in L^1(G)$ and for any $m,n\in
\Z^\sigma$, $\what{f}(\sigma, \lambda)_{m,n}$ can be defined for
any $\lambda\in \C$. We will use the following asymptotic
behaviour of $\phi_\lambda$ (see \cite[p. 447]{He2}):
\begin{equation}
\lim_{t\to \infty}
e^{(-|\lambda|+1)t}\phi_{|\lambda|}(a_t)=c(|\lambda|)\mbox{ for }
\lambda\neq 0 \label{asymp-phi}
\end{equation}
 where
$c(\lambda)=\frac{\Gamma(\lambda/2)}{\Gamma(1+\lambda/2)}$
(\cite[p. ~24]{Ba}) is the Harish-Chandra $c$-function for $G$.

Let $B\subset G$ be a large compact set containing the identity.
Since $f\in L^2(G)$, $f$ is a locally integrable function on $G$
and hence $\int_B|f(x)|dx<\infty$ and
$\int_B|f(x)|\Phi^{m,n}_{\sigma,\lambda}(x)dx<\infty$.

We claim that $\what{f}(\sigma, \cdot)$ cannot be supported on a
set of finite measure in $i\R$. Suppose $\what{f}(\sigma, \cdot)$
is supported on a set of finite measure. Then  for any $m,n\in
\Z^\sigma$ $\what{f}(\sigma, \lambda)_{m,n}$ is also supported on
a set of finite measure.  Now as $\what{f}(\sigma,
\cdot)\not\equiv 0$, there exists $\lambda_0\neq 0$ such that
$\what{f}(\sigma, \lambda_0)\neq 0$ and from
(\ref{beurling-condition-0-main})
$\int_G|f(x)|\phi_{|\lambda_0|}(x)dx<\infty$. Suppose
$|\lambda_0|=r>0$. Using (\ref{asymp-phi}) we have
$\int_{G\setminus B}|f(x)|e^{(r-1)\sigma(x)}dx<\infty$. Again
using the fact that $f$ is locally integrable, we see that
$\int_B|f(x)|e^{(r-1)\sigma(x)}dx<\infty$, for  $r$ as above.

Together they imply, $\int_G|f(x)|e^{(r-1)\sigma(x)}dx<\infty$.
Using the Cartan decomposition $G=K\overline{A^+}K$, we can
rewrite it as
\begin{equation}
\int_{K}\int_0^\infty\int_K|f(k_1a_tk_2)|e^{(r-1)t}\sinh
2t\,dk_1\,dt\,dk_2 <\infty. \label{has-decay-1}
\end{equation}

We have $\what{f}_{m,n}(\sigma, \lambda)=\int_Gf(x)
\Phi^{m,n}_{\lambda}(x^{-1})dx$. Then
$$\ba{lll}|\int_Gf(x) \Phi^{m,n}_{\lambda}(x^{-1})dx|&\le& \int_G|f(x)||\Phi^{m,n}_{\lambda}(x^{-1})|dx\\ \\
&\le&\int_G|f(x)|e^{|\Re \lambda|\sigma(x)}\Xi(x)dx.\\
\\ &\le&\int_{K}\int_0^\infty\int_K|f(k_1a_tk_2)|e^{(|\Re
\lambda|-1)t}(1+t)\sinh 2t\, dk_1\,dt\,dk_2.\ea$$ In the last two
steps we have used the estimates of $\Phi^{m,n}_{\lambda}(x)$ and
$\Xi(x)$ (see (\ref{basic-estimates})) and the Cartan
decomposition. Thus we have for $0<r'<r$ and for $\lambda\in \C$
with $\Re \lambda|<r'$,
$$|\int_G f(x)\Phi^{m,n}_{\lambda}(x^{-1})dx|\le
\int_{K}\int_0^\infty\int_K|f(k_1a_tk_2)|e^{(r-1)t}e^{(|\Re
\lambda|-r')t}e^{(r'-r)t}(1+t)\sinh 2t\, dk_1\,dt\,dk_2.$$
Hence
by (\ref{has-decay-1}), $|\int_G
f(x)\Phi^{m,n}_{\lambda}(x^{-1})dx|<\infty$ for any $\lambda\in
\C$ with $|\Re \lambda|<r'$.

By a standard use of Morera's theorem it follows that
$\what{f}_{m,n}$ is analytic in the open strip $|\Re \lambda|<r'$
in $\mathfrak a^*_\C=\C$. This
 contradicts the assumption that $\what{f}(\sigma, \cdot)$ and hence $\what{f}_{m,n}(\sigma, \cdot)$ is
 supported on a set of finite measure. Thus our claim is established, i.e.
 $\what{f}(\sigma, \cdot)$ is supported on a set of infinite measure.

Now as $\what{f}(\sigma, \cdot)$ is supported on a set of infinite
measure, from (\ref{beurling-condition-0-main}) and
(\ref{asymp-phi}), it follows that for any large $M>0$ there
exists $\lambda\in i\R$, $|\lambda|>M$ such that $\int_{G\setminus
B} |f(x)|e^{(|\lambda|-1)\sigma(x)}dx<\infty$. This implies that
$\int_{G\setminus B}|f(x)|dx<\infty$. Since
$\int_B|f(x)|dx<\infty$ we immediately see that $f\in L^1(G)$.

Now given any $\lambda'\in \C$ with $|\lambda'|=l$ say, we choose
$M=l$ in the above. Then $\int_{G\setminus
B}|f(x)|e^{(|\lambda|-1)\sigma(x)}dx<\infty$ for some $\lambda$
with $|\lambda|>l$. Since (see section 2)
$$|\Phi^{m,n}_{\sigma, \lambda'}(x)|\le \phi_{\Re \lambda'}(x)\le e^{|\Re
\lambda'|\sigma(x)}\Xi(x)\le e^{(|\Re
\lambda'|-1)\sigma(x)}(1+\sigma(x)),\mbox { for any }\lambda'\in
\C$$ we see that $\int_{G\setminus B} |f(x)||\Phi^{m,n}_{\sigma,
\lambda'}(x)|dx<\infty$. Combining this with the fact
$\int_B|f(x)||\Phi^{m,n}_{\sigma, \lambda'}(x)|dx<\infty$, we have
$\int_G|f(x)||\Phi^{m,n}_{\lambda'}(x)|dx<\infty$. Thus we have
established that the Fourier transform $\what{f}(\sigma,
\lambda)_{m,n}$ of $f$ exists for every $\lambda\in \C$.  Notice
that we have actually established that the function $|f|$ has
enough decay so that for any $m,n\in \Z^\sigma$,
$|f|\,\what{\,}\,(\sigma, \lambda)_{m,n}$ exists for any
$\lambda\in \C$. \vspace{.2in}

\noindent{\bf Step 2:} We have  $f=\sum_{m,n\in \Z}f_{m,n}$ in the
sense of distributions on $G$  (see section 2). Note that for each
$m,n\in \Z$, $f_{m,n}$ is in $L^1(G)\cap L^2(G)$, since $f\in
L^1(G)\cap L^2(G)$. Let us fix a $\sigma\in \what{M}$ and take two
arbitrary $m,n\in \Z^\sigma$. As $|f_{m,n}(x)|\le \int_{K\times
K}|f(k_1xk_2)|dk_1\,dk_2$, $\phi_{|\lambda|}(x)$ is
$K$-biinvariant and the Haar measure $dx$ is invariant under the
transformation $x\mapsto k_1xk_2$, we can substitute $f$ by
$f_{m,n}$ in (\ref{beurling-condition-0-main}) and in
(\ref{beurling-discrete-condition-main}). Also in
(\ref{beurling-condition-0-main}) (respectively
(\ref{beurling-discrete-condition-main})) we can substitute
$\|\what{f}(\sigma, \lambda)\|_2$ by $|\what{f}_{m,n}(\lambda)|$
(respectively $\|\what{f}(k)\|_2$ by $|\what{f}_{m,n}(k)|$) as
$\|\what{f}(\sigma,
\lambda)\|_2^2=\sum_{m,n}|\what{f}_{m,n}(\sigma,\lambda)|^2$
(respectively
$\|\what{f}(k)\|_2^2=\sum_{m,n}|\what{f}_{m,n}(k)|^2$ ). Thus we
get
\begin{equation}\ba{l}\int_G\int_{i\R}|f_{m,n}(x)||\what{f}_{m,n}(\sigma,
\lambda)|\phi_{|\lambda|}(x)dxd\mu(\sigma,\lambda)< \infty,\\
\\
\sum_{k\in\Z^*}\frac{|k|}{2\pi}\int_G|f_{m,n}(x)||\what{f}_{m,n}(k)|\phi_{|k|}(x)dx<\infty.
\ea \label{beurling-condition-0-mn}
\end{equation}

Starting from (\ref{beurling-condition-0-mn}) if we can show that
$f_{m,n}=0$ then we are done in view of the decomposition of $f$
in $f_{m,n}$. So, we can confine ourselves to the set of functions
of type $(m,n)$ for some $m,n\in \Z$ of the same parity.

In order to avoid complicated notation we will simply write $\ff$
for $f_{m,n}$. Also by $\what{\ff}(\lambda)$ (respectively
$\what{\ff}(k)$) we will mean $\what{f}(\lambda)_{m,n}$
(respectively $\what{f}(k)_{m,n}$).
Notice that we have omitted $\sigma$ and have written
$\what{\ff}(\lambda)$ for $\what{\ff}(\sigma, \lambda)$ as the
$\sigma\in \what{M}$  is fixed when $m,n$ are fixed by $m,n\in
\Z^\sigma$. For the same reason we will write $\mu(\lambda)$ for
$\mu(\sigma, \lambda)$. So we rewrite the inequalities
(\ref{beurling-condition-0-mn}) as:
\begin{equation}\ba{l}\int_G\int_{i\R}|\ff(x)||\what{\ff}(
\lambda)|\phi_{|\lambda|}(x)dxd\mu(\lambda)< \infty,\\ \\
\sum_{k\in\Z^*}\frac{|k|}{2\pi}\int_G|\ff(x)||\what{\ff}(k)|\phi_{|k|}(x)dx<\infty,
\ea \label{beurling-condition-0-mn-2}
\end{equation}
where $\ff\in L^1(G)\cap L^2(G)$ and $\ff$ is of type $(m,n)$.

As
$|\ff|\,\what{\,}\,(|\lambda|)=\int_G|\ff(x)|\phi_{|\lambda|}(x)dx$,
the  inequality above is equivalent to
$$\int_{i\R}|\ff|\,\what{\,}\,(|\lambda|)\,|\what{\ff}(\lambda)|d\mu(\lambda)
=\int_\R
|\ff|\,\what{\,}\,(|\lambda|)\,|\what{\ff}(i\lambda)|d\mu(i\lambda)<\infty.$$
(Note that $|\ff|\,\what{\,}\,(|\lambda|)$ exists for $\lambda\in
i\R$ for reasons mentioned in Step 1.)

Using lemma 3.4 we get
$$\int_\R\wtilde{{\mathcal A}
|\ff|}(-i|\lambda|)\,|\wtilde{{\mathcal A}
\ff}(\lambda)|d\mu(i\lambda)<\infty.$$ Recall that $\ff$ being a
$K$-finite function of type $(m,n)$, $|\ff|$ is $K$-biinvariant
and hence ${\mathcal A}|\ff|$ is an even function. Therefore from
above we have
\begin{equation}\int_\R\wtilde{{\mathcal A}
|\ff|}(i|\lambda|)\,|\wtilde{{\mathcal A}
\ff}(\lambda)|d\mu(i\lambda)<\infty. \label{step1}
\end{equation}

 Now $\wtilde{{\mathcal A}
|\ff|}(i|\lambda|)=\int_\R{\mathcal A}|\ff|(t)e^{|\lambda|t}dt\ge
\int_0^\infty{\mathcal A}|\ff|(t)e^{|\lambda|t}dt=\frac
12M({\mathcal A} |\ff|)(\lambda)$, since $A|\ff|$ is an even
function. Here $M({\mathcal A} |\ff|)(\lambda)$ is as defined in
(\ref{modified-hormander}).

 So from (\ref{step1}) we have $$\int_\R M({\mathcal A}
|\ff|)(\lambda)\,|\wtilde{{\mathcal A}
\ff}(\lambda)|d\mu(i\lambda)<\infty$$ As $|{\mathcal A} \ff(t)|\le
{\mathcal A}|\ff|(t)$ for all $t\in\R$ (see (\ref{mod-f-abel})),
we have $M({\mathcal A} \ff)(\lambda)\le M({\mathcal A}
|\ff|)(\lambda)$ since $e^{|\lambda||t|}$ is positive.

Therefore \begin{equation}\int_\R M({\mathcal A}
\ff)(\lambda)\,|\wtilde{{\mathcal A}
\ff}(\lambda)|d\mu(i\lambda)<\infty.\label{step2}\end{equation}
\vspace{.2in}

\noindent{\bf Step 3} In this step we will first show that in
(\ref{step2}) the Plancherel measure $d\mu(i\lambda)$ can be
substituted  by the Lebesgue measure $d\lambda$ and then conclude
that $\what{\ff}(\lambda)=0$ for all $\lambda\in i\R$.

We have  $\int_{\R}M({\mathcal A}\ff)(\lambda) |\wtilde{{\mathcal
A }\ff}(\lambda)|\mu(i\lambda) d\lambda <\infty$ since
$d\mu(i\lambda)=\mu(i\lambda)d\lambda$ where  $\mu(i\lambda)$ is
as in (\ref{plancherelmeasure}). Therefore $M({\mathcal
A}\ff)(\lambda)|\wtilde{{\mathcal A}\ff}(\lambda)| \mu(i\lambda)$
is finite for almost every $\lambda\in \R$ with respect to the
Lebesgue measure. Now $\mu(i\lambda)$ and
$\what{\ff}(i\lambda)=\wtilde{{\mathcal A}\ff}(\lambda)$ are real
analytic functions for $\lambda\in \R$. Hence $M({\mathcal
A}\ff)(\lambda)$ is finite almost everywhere. Now by its very
definition $M({\mathcal A}\ff)(\lambda)$ is even in $\lambda$ and
an increasing function of $|\lambda|$. Consequently $M({\mathcal
A}\ff)(\lambda)$ is finite everywhere and locally integrable,
because for any $R^{\prime}>0 $ we have $\int_{|\lambda|\leq
R^{\prime}} M({\mathcal A}\ff)(\lambda) d\lambda \leq 2 R^{\prime}
M(R^{\prime})< \infty$. From (\ref{plancherelmeasure}) we see that
there exists $R > 0 $ such that $\mu(i\lambda)$ is $\geq \f
{1}{2}$ say for all $\lambda\in \R,|\lambda| \geq R $ because $
\mu(i\lambda) \rightarrow {\iy}$ as $|\lambda|\rightarrow {\iy}$.
The Euclidean Riemann-Lebesgue lemma shows that $\wtilde{\mathcal
A \ff}(\lambda)$ is a bounded continuous function. Therefore
$\int_{\R} M({\mathcal A}\ff)(\lambda)|\wtilde{\mathcal A
\ff}(\lambda)|d\lambda$ is finite. This shows ${\mathcal A}\ff$
satisfies the condition (\ref{modified-hormander}). Consequently
$\wtilde{{\mathcal A}\ff}(\lambda)=\what{\ff}(i\lambda)= 0$ for
all $\lambda \in \R $. That is $\what{\ff}(\lambda)=0$ for all
$\lambda\in i\R$.  Hence the Fourier transform of $\ff$ is
supported on the discrete series. \vspace{.2in}

\noindent{\bf Step 4} Let $D_{m,n}=\{k \in \Z^* |\ \pi_k \mbox{
has } m \mbox{ and } n\mbox{ as } K \mbox{ types}\}$. That is
$D_{m,n}$ is the set parametrizing the discrete series
representations which admit the pair $(m,n)$ as K types. Note that
the cardinality of $D_{m,n}$ is finite. Since $\ff$ is of type
$(m,n)$ and its Fourier transform is supported only on the
discrete series, it must be a finite linear combination of matrix
coefficients of discrete series representations parametrized by
elements of $D_{m,n}$. That is $\ff=\Si_{k \in D_{m,n}}
c_k\Psi^{m,n}_k$, where $\Psi^{m,n}_k$ is the $(m,n)$-th canonical
matrix coefficient of $\pi_k$. Let $k_0=\max \{|k|\ |\ k\in
D_{m,n}\}$.  Now (see \cite[p. ~70]{Ba}) $\Psi^{m,n}_k(a_t)$ is
asymptotic to $A e^{-(1 + |k|)t}$ as $t\longrightarrow\infty$
where $A$ is a nonzero constant depending on $k$. This shows that
$|\ff(a_t)|\geq A' e^{-(1+k_0)t}$ for all $t>0$ sufficiently
large, where $A'$ is a positive constant depending on $\ff$. Now
using Cartan decomposition and (\ref{asymp-phi}) we see that
$\int_G |\ff(x)|\phi_{k_0}(x) dx=\int_0^\infty
|\ff(a_t)|\phi_{k_0}(a_t)\sinh 2t dt$ is infinite, where $k_0$ is
as above. Therefore $\ff$ cannot satisfy the second condition of
the theorem unless it is zero almost everywhere.

Thus we have shown that $\ff(x)=f_{m,n}(x)=0$ for almost all $x\in
G$. As $(m,n)$ is arbitrary, in view of the decomposition
$f=\sum_{m,n}f_{m,n}$ in the sense of distributions, it follows
that $f(x)=0$ for almost all $x\in G$.
\end{proof}
\section{Consequences}\newsection
We have already mentioned that our theorem \ref{thm-sl2-1} implies
the other QUP's. We now state these  (see for instance \cite{SiSu,
Sa1,Se, NR, FS, T1} for independent proofs of these theorems).

1. Hardy's Theorem: Let $f:G \rightarrow \C$ be a complex valued
measurable function and assume  that,
$$\ba{lll} (1) & |f(x)| \leq Ce^{-\alpha \s(x)^2}& \mbox{ for all } x\in G,\\ &\\
(2) & \|\what{f}{(\s, \lambda)}\|_2 \leq Ce^{-\beta|\lambda|^2}&
\mbox{ for all } \sigma\in \what{M}\mbox{ and } \lambda\in
i\R,\ea$$ where $\alpha, \beta$ are positive constants. If $\alpha
\beta > \f{1}{4} $ then $f=0$ almost everywhere.

2. Morgan's Theorem (strong version): Let $f:G \rightarrow \C$ be
measurable and assume that,
$$\ba{lll} (1)& |f(x)| \leq Ce^{-\alpha \s(x)^p}& \mbox{ for all } x\in G,\\&\\
(2)& \|\what{f}{(\sigma, \lambda)}\|_2 \leq e^{-\beta|\lambda|^q}
& \mbox{ for all } \sigma\in \what{M}\mbox{ and } \lambda\in
i\R,\ea$$ where $\alpha, \beta $ are positive constants, $1< p <
{\iy}$ and
$\f 1{p} + \f{1}{q} = 1 $.\\
If $(\alpha p )^{\f{1}{p}}(\beta q)^{\f{1}{q}} >1 $, then $f= 0$
almost everywhere.

3. Cowling-Price Theorem: Let $f: G \rightarrow \C$ be measurable
and assume that for positive constants $\alpha$ and $\beta$ we
have
$$\ba{ll} (1)& e_{\alpha}f \in L^p(G),\\&\\
(2)& e_{\beta}\|\what{f}(\s, \lambda)\|_2 \in L^q(i\R;
d\mu(\sigma, \lambda)),\ea$$ where $ 1 \leq p,q \leq{\iy}$,
$e_{\alpha}(x) = e^{\alpha{\sigma}^2(x)}$ and $
e_{\beta}(\lambda)= e^{\beta |\lambda|^2}$.

If $\alpha\beta>\frac 14$, then $f=0$ almost everywhere.

4. Gelfand-Shilov: Let $f\in L^2(G)$. Suppose $f$ satisfies
$$\ba{ll} (1)& \int_G|f(x)|e^{\frac{(\alpha\sigma(x))^p}{p}}dx<\infty,\\&\\
(2)& \int_{i\R}\|\what{f}(\sigma,\lambda)\|_2
e^{\frac{(\beta|\lambda|)^q}{q}}d\mu(\sigma, \lambda)<\infty,\ea$$
where $ 1 < p,q <{\infty}$, $\frac 1p+\frac 1q=1$ and
$\alpha\beta\ge 1$. Then $f=0$ almost everywhere.

The proof of the deduction of these theorems from theorem
\ref{thm-sl2-1} is similar to that in the Euclidean case. We
illustrate this in the case of Cowling-Price theorem. As in the
proof of theorem \ref{thm-sl2-1} we can reduce the proof to the
case when $f$ is of type $(m,n)$. We will write
$\what{f}(\lambda)$ for $\what{f}(\sigma, \lambda)$ and
$\mu(\lambda)$ for $\mu(\sigma, \lambda)$ as  $\sigma\in \what{M}$
is uniquely determined by $m,n\in \Z^\sigma$.

Let $f$ be a function of type $(m,n)$ which satisfies the
conditions (1) and (2) of  the Cowling-Price theorem. Then we can
choose $0<\alpha'<\alpha$ (respectively $0<\beta'<\beta$) such
that $\alpha'\beta'>\frac 14$. We have $e_{\alpha'}|f|\in L^1(G)$
(respectively $e_{\beta'}|\what{f}|\in L^1(i\R, d\mu(\lambda)$).
We will  show that
\begin{equation}\int_G\int_{i\R} |f(x)| |\what{f}(\lambda)|
e^{\sigma(x)|\lambda|} dx d \mu(\lambda)<\infty.\label{C-P}
\end{equation}

We take $\beta''<\beta'$ such that $\alpha'\beta''=\frac 14$. Then
$e_{\beta''}\what{f}\in L^1(i\R, d\mu(\lambda))$ and
$$\ba{lll}\tilde{I}&=&\int_G\int_{i\R}e_{\alpha'}(x)|f(x)|\,e_{\beta''}(\lambda)|\what{f}(\lambda)|
e^{-\alpha'\sigma(x)^2}\,e^{-\beta''|\lambda|^2}e^{\sigma(x)|\lambda|}dxd\mu(\lambda)\\
&&\\
&=&\int_G\int_{i\R}e_{\alpha'}(x)|f(x)|e_{\beta''}(\lambda)|\what{f}(\lambda)|
e^{-(\sqrt{\alpha'}\sigma(x)-\sqrt{\beta''}|\lambda|)^2}dxd\mu(\lambda).
\ea$$ Since
$e^{-(\sqrt{\alpha'}\sigma(x)-\sqrt{\beta''}|\lambda|)^2}\le 1$,
$\tilde{I}<\infty$.

 Furthermore the rapid decay of $f$ namely $e_{\alpha'}f(x)\in
L^1(G)$ immediately shows that $f$ also satisfies the second
condition of theorem \ref{thm-sl2-1}. That is the pair $(f,
\what{f})$ satisfies the condition of theorem \ref{thm-sl2-1}.
Hence $f=0$ almost everywhere.

In view of the inequality $\alpha\beta\sigma(x)|\lambda|\le
\frac{\alpha^p}{p}\sigma(x)^p+\frac{\beta^q}{q}|\lambda|^q$
conditions (1) and (2) of the Gelfand-Shilov theorem immediately
imply the first condition of theorem \ref{thm-sl2-1}. Again the
rapid decay of $f$ given in condition (1) of the Gelfand-Shilov
theorem shows that $f$ also satisfies the second condition of
theorem \ref{thm-sl2-1}.

Morgan's theorem easily follows from Gelfand-Shilov theorem. Note
that the special case $ p=q=2$ in Morgan's theorem gives Hardy's
theorem.

\begin{Remark}
{\em A careful reader will observe that this is the first time
when in a theorem of uncertainty, discrete series representations
appear in the hypothesis; compare with for instance \cite{CSS,
Sa1}. It is obvious that for an integrable cusp form $\Psi$, its
Fourier transform vanishes identically on the unitary principal
series and hence $\Psi$ trivially satisfies the first condition of
Beurling's theorem. The non-appearance of the discrete series in
the other QUP's mentioned in  this section can be explained by
noting that in all of them we put very rapid decay on the function
$f$ which forces every $K$-finite component of $f$ to satisfy the
second condition of the Beurling's theorem (see
(\ref{beurling-discrete-condition-main}) and
(\ref{beurling-discrete-condition-main-reform1})).

Our result also indicates that for a group having real rank
greater than 1, the hypothesis of Beurling's theorem will involve
all non-minimal principal series and discrete series, in contrast
with the other QUP's.}\end{Remark}

\end{document}